\documentclass{article}
\usepackage{amsmath}
\usepackage{amssymb}

\newcommand{\TT}{\mathbf{T}}
\newcommand{\UU}{\mathbf{U}}
\newcommand{\kk}{\mathbf{k}}
\newcommand{\zz}{\mathbf{z}}
\newcommand{\WW}{\mathbf{W}}

\newcommand{\T}{\mathbb{T}}

\newcommand{\C}{\mathbb{C}}
\newcommand{\R}{\mathbb{R}}
\newcommand{\Z}{\mathbb{Z}}
\newcommand{\N}{\mathbb{N}}

\def\co{\colon}
\def\ba{^{\relbar}}
\def\la{\langle}
\def\ra{\rangle}
\def\i{\infty}

\def\op{\oplus}

\def\ol{\overline}

\def\HH{{\cal H}}
\def\L{{\cal L}}
\def\M{{\cal M}}
\def\K{{\cal K}}
\def\I{{\cal I}}
\def\E{{\cal E}}
\def\F{{\cal F}}

\def\P{{\cal P}}

\def\B{{\cal B}}
\def\NN{{\cal N}}

\def\b{\beta}
\def\g{\gamma}
\def\f{\varphi}

\def\o{\omega}
\def\O{\Omega}

\def\z{\zeta}

\def\s{\sigma}

\def\th{\vartheta}
\def\FF{\Phi}
\def\PP{\Psi}
\def\SS{\Sigma}
\def\h{\chi}
\def\d{\delta}
\def\l{\lambda}
\def\m{\mu}
\def\k{\kappa}

\def\eq{\simeq}
\def\ae{\overset{a}\sim}
\def\ac{\overset{a}\prec}
\def\wt{\widetilde}
\def\oinf{\textup{o-}\inf}
\def\Mlim{M\textup{-}\lim}
\def\Hlat{\hbox{ Hlat }}
\def\eqE{\overset{E}=}
\def\subE{\overset{E}\subset}
\def\veeE{\overset{E}\vee}

\newtheorem{pro}{Proposition}
\newtheorem{thm}[pro]{Theorem}
\newtheorem{lem}[pro]{Lemma}
\newtheorem{cor}[pro]{Corollary}
\newtheorem{ex}[pro]{Example}
\newtheorem{rem}[pro]{Remark}
\newtheorem{que}[pro]{Question}

\begin{document}

\title{Quasianalytic $n$-tuples of Hilbert space operators}

\author{L\'aszl\'o K\'erchy}

\date{}

\maketitle

\begin{abstract}
\noindent
In this paper a systematic study of unitary asymptotes of commuting $n$-tuples of general Hilbert space operators is initiated.
Special emphasis is put on the study of the quasianalicity property.


\end{abstract}

\section{Introduction} 
\label{intr}

The theory of Hilbert space operators contains many beautiful results, yet some of its fundamental questions --- like the hyperinvariant subspace problem --- are still open; see, e.g., the monographs \cite{ConwayFA}, \cite{ConwayOT}, \cite{Radjavi} and \cite{Chalendar}.
Beside the structure of single operators, multivariable operator theory also came into the focus of researchers' interest; see, e.g., \cite{Douglas}, \cite{Chen} and \cite{Muller}.
The spectral analysis of unitary operators makes possible to explore their structure.
Hence, it has been reasonable to relate more general operators to unitaries.
Such connection has been made for contractions by Sz.-Nagy's dilation theorem; the resulting theory is presented in \cite{NFBK}.
The residual and $*$-residual parts of the unitary dilation proved to be especially useful in the study of contractions.
A more direct approach to these components, originated in \cite{Nagy47}, leads to the concept of unitary asymptote, and opens the way for generalizations to more general settings.
The fundamental properties of unitary asymptotes have been summarized for contractions in Chapter IX of \cite{NFBK}, for power bounded operators in \cite{KerInd89}, for commuting $n$-tuples of power bounded operators in \cite{Berc}, for representations of abelian semigroups in \cite{KerDiscrepr} and \cite{KerLeka}; see also the references therein.
In the present paper we initiate a systematic study of unitary asymptotes for commuting $n$-tuples of general Hilbert space operators, following the categorical approach applied in \cite{BK} and Chapter IX of \cite{NFBK}.
In the quest for proper hyperinvariant subspaces of asymptotically non-vanishing contractions, the crucial class is formed by the quasianalytic contractions, which show homogeneous behaviour in localization.
For their study see, e.g., \cite{Ker2001}, \cite{Ker2011}, \cite{KerTot}, \cite{KSzStudia} and \cite{KSzProc}.
These investigations were extended to polynomially bounded operators in \cite{KerOT} and \cite{KerActa}; see also \cite{Gamal}.
In this paper we provide a systematic study of quasianalycity in the setting of commuting $n$-tuples of general Hilbert space operators.

Our paper is organized in the following way.
In Section \ref{unitasymp} the concepts of unitary intertwining pairs and unitary asymptotes are introduced.
Fundamental properties, like minimality, norm-control and uniqueness are examined.
Necessary condition is given for the existence of the unitary asymptote in terms of the spectral radius.
Adjoints and restrictions to invariant subspaces are considered.
In Section \ref{orbits} it is shown that the Lower Orbit Condition (LOC) is a sufficient but not necessary condition of the existence of unitary asymptotes.
Application of invariant means yields that (LOC) holds in the power bounded case.
A classification is made relying on the annihilating subspace.
Section \ref{quasianalytic} deals with quasianalycity and hyperinvariant subspaces.
The commutant mapping relates the commutant of the $n$-tuple $\TT=(T_1,\dots,T_n)$ to the commutant of the corresponding $n$-tuple $\UU=(U_1,\dots,U_n)$ of unitaries.
We obtain that hyperinvariant subspaces of $\UU$ induce hyperinvariant subspaces of $\TT$.
The basic facts of the spectral analysis of $\UU$ are recalled; in particular, the spectral subspaces are hyperinvariant.
Quasianalycity of $\TT$ means that the localizations of the spectral measure $E$ of $\UU$ at non-zero vectors, corresponding to vectors in the space of $\TT$, are equivalent.
If this homogeneity breaks, then $\TT$ has proper hyperinvariant subspaces.
In Section \ref{sp-invariants} the local residual sets are introduced.
Exploiting the lattice structure of the Borel sets on $\T^n$, the quasianalytic spectral set $\pi(\TT)$ is defined.
Quasianalycity is related to a cyclicity property of the commutant of $\UU$.
The absolutely continuous (a.c.) case is also studied.
The a.c.\ global residual set $\o_a(\TT)$ is the measurable support of the spectral measure $E$, and quasianalycity is equivalent to the coincidence $\pi_a(\TT)= \o_a(\TT)$.

We shall use the following notation: $\N, \Z_+, \Z, \R_+, \R$ and $\C$ stand for the set of positive integers, non-negative integers, integers, non-negative real numbers, real numbers and complex numbers, respectively.
For any $n\in\N, \; \N_n:= \{k\in\N: k\le n\}$, and $\T:= \{z\in\C: |z|=1\}$.

\section{Existence of unitary asymptotes}
\label{unitasymp}

Let $\TT=(T_1,\dots,T_n)$ be a commuting $n$-tuple of (bounded, linear) operators on the (complex, separable) Hilbert space $\HH \; (n\in\N)$.
Giving any other commuting $n$-tuple $\wt\TT= (\wt T_1,\dots,\wt T_n)$ of operators on the Hilbert space $\wt\HH$, the \emph{intertwining set} $\I(\TT,\wt\TT)$ consists of those (bounded, linear) transformations $X\co \HH\to\wt\HH$, which satisfy the conditions $XT_i=\wt T_iX\; (i\in\N_n)$.
For any $X\in \I(\TT,\wt\TT)$ and $\kk=(k_1,\dots,k_n)\in \Z_+^n, \; X\TT^\kk = \wt\TT^\kk X$ holds, where $\TT^\kk:= T_1^{k_1} \cdots T_n^{k_n}$.
Hence $Xp(\TT)= p(\wt\TT)X$ is true for every polynomial $p(\zz)= \sum_\kk c_\kk \zz^\kk$ in the variable $\zz=(z_1,\dots,z_n)\in \C^n$, where $\zz^\kk:= z_1^{k_1}\cdots z_n^{k_n}$.
It is clear that $\I(\TT,\wt\TT)$ is a Banach space with the operator norm.
The \emph{commutant} $\{\TT\}':= \I(\TT,\TT)$ of $\TT$ is a Banach subalgebra of the operator algebra $\L(\HH)$.

Let $\UU=(U_1,\dots,U_n)$ be a commuting $n$-tuple of unitary operators on the Hilbert space $\K$.
If $X\in \I(\TT,\UU)$, then $(X,\UU)$ is called a \emph{unitary intertwining pair} of $\TT$.
The subspace $\K_0= \vee\{\UU^\kk X\HH: \kk\in\Z^n\}$ is the smallest one, reducing $\UU$ and containing the range of $X$.
The restriction $\UU_0= (U_{10}, \dots, U_{n0})$, where $U_{i0}= U_i|\K_0\; (i\in\N_n)$, is a commuting $n$-tuple of unitaries on $\K_0$; we use the notation: $\UU_0= \UU|\K_0$.
Setting $X_0\co \HH\to \K_0, h\mapsto Xh$, the unitary intertwining pair $(X_0,\UU_0)$ of $\TT$ is called the \emph{minimal part} of $(X,\UU)$, and the notation $(X,\UU)_0:= (X_0,\UU_0)$ is used.
We say that $(X,\UU)$ is \emph{minimal}, if $\K_0=\K$.

The unitary intertwining pairs $(X,\UU)$ and $(\wt X, \wt\UU)$ of $\TT$ are called \emph{similar}, in notation: $(X,\UU) \approx (\wt X,\wt\UU)$, if there exists an invertible $Z\in \I(\UU,\wt\UU)$ satisfying the condition $\wt X= ZX$.
These pairs are called \emph{equivalent}, in notation: $(X,\UU) \eq (\wt X, \wt\UU)$, if there exists a unitary $Z\in \I(\UU,\wt\UU)$ such that $\wt X=ZX$.

The unitary intertwining pair $(X,\UU)$ of $\TT$ is a \emph{unitary asymptote} of $\TT$, if it is universal in the sense that for any other unitary intertwining pair $(X', \UU')$ of $\TT$ there exists a unique $Y'\in \I(\UU,\UU')$ such that $X'=Y'X$.

\begin{pro}
\label{minimal}
If $(X,\UU)$ is a unitary asymptote of $\TT$, then $(X,\UU)$ is minimal.
\end{pro}

\noindent
\textbf{Proof.}
Set $\K_0=\vee\{\UU^\kk X\HH: \kk\in \Z^n\}$ as before.
Since $IX=X= P_{\K_0}X$ and $I, P_{\K_0}\in \I(\UU,\UU)$, it follows by the uniqueness condition that $I=P_{\K_0}$, and so $\K_0=\K$.

\rightline{$\square$}

If $\UU'$ is a commuting $n$-tuple of unitaries on $\K'$ and $Y'\in \I(\UU,\UU')$, then $(X'=Y'X,\UU')$ is obviously a unitary intertwining pair of $\TT$ and $\|X'\|\le \|Y'\|\cdot\|X\|$.
We say that the unitary asymptote $(X,\UU)$ of $\TT$ has \emph{norm-control} with $\k\in\R_+$, if $\|Y'\|\le \k\|Y'X\|$ holds for every unitary intertwining pair $(Y',\UU')$ of $\UU$.
The smallest possible $\k$ is called the \emph{optimal norm-control}, and it is denoted by $\k_{\rm op}=\k_{\rm op}(X,\UU)$.

The following proposition shows that the minimal unitary intertwining pairs of $\TT$ can be obtained, up to equivalence, by the aid of operators in the commutant of $\UU$.

\begin{pro}
\label{generating}
Assume that $(X,\UU)$ is a unitary asymptote of $\TT, \;\UU$ acting on $\K$.
Let $(X',\UU')$ be a unitary intertwining pair of $\TT$, let $(X',\UU')_0= (X'_0,\UU'_0)$ be its minimal part, where $X'=Y'X$ and $X'_0=Y'_0X$ with unique $Y'\in\I(\UU,\UU')$ and $Y'_0=\I(\UU,\UU'_0)$.
Then
\begin{itemize}
\item[\textup{(a)}] $\; Y'k=Y'_0k\;$ for all $k\in\K$;
\item[\textup{(b)}] $\; (X',\UU')_0\eq (|Y'|X,\UU)_0$, where $|Y'|\in\{\UU\}'$.
\end{itemize}
\end{pro}

\noindent
\textbf{Proof.} 
The $n$-tuples $\UU'$ and $\UU'_0$ act on $\K'$ and $\K'_0$, respectively.
For every $h\in\HH$, we have $Y'Xh=X'h\in\K'_0$. 
Since $\K'_0$ is reducing for $\UU'$, it follows that $Y'\UU^\kk Xh= \UU'^\kk Y'Xh\in \K'_0$ for all $\kk\in\Z^n$.
By the minimality of $(X,\UU)$ we infer that $Y'\K\subset \K'_0$.
Setting $Y'_{00}\co \K\to\K'_0, k\mapsto Y'k$, it is clear that $Y'_{00}\in \I(\UU,\UU'_0)$.
Furthermore, the equalities $Y'_{00}Xh =Y'Xh= X'h=X'_0h= Y'_0Xh\; (h\in\HH)$ yield that $Y'_{00}=Y'_0$, and so (a) is proved.

In view of (a) we know that $|Y'|= |Y'_0|$.
Let us consider the polar decomposition $Y'_0= W'_0|Y'_0|$, where $W'_0$ is a partial isometry with $\ker W'_0= \ker|Y'_0|$.
Since $Y'_0\in \I(\UU,\UU'_0)$ readily implies 
$Y'_0\in\I(\UU^*,{\UU'_0}^*)$, 
it can be easily seen that $|Y'_0|\in \{\UU\}'$ and $W'_0\in \I(\UU, \UU'_0)$.
Furthermore, the subspaces $\K^{(0)}= (\ker W'_0)^\perp$ and $W'_0\K^{(0)}$ are reducing for $\UU$ and $\UU'_0$, respectively.
Since $W'_0\K^{(0)}\supset W'_0|Y'_0|\K= Y'_0\K\supset Y'_0X\HH= X'_0\HH$, it follows that $W'_0\K^{(0)} =\K'_0$.
It is clear that $W^{(0)}\in \I(\UU^{(0)}, \UU'_0)$ is a unitary transformation, where $\UU^{(0)}:= \UU|\K^{(0)}$ and $W^{(0)}:= W'_0|\K^{(0)}$.
On the other hand, the minimal part of the unitary intertwining pair $(|Y'_0|X, \UU)$ of $\TT$ is $(X^{(0)},\UU^{(0)})$, where $X^{(0)}\co \HH\to \K^{(0)}, h\mapsto |Y'_0|Xh$.
For every $h\in\HH$, we have $W^{(0)} X^{(0)}h = W'_0|Y'_0|Xh = Y'_0Xh= X'_0h$; hence $W^{(0)} X^{(0)}= X'_0$.
Therefore, $(X^{(0)}, \UU^{(0)})$ is equivalent to $(X'_0,\UU'_0)$.

\rightline{$\square$}

We show that having a norm-control is a general property.

\begin{pro}
\label{control}
Assume that $(X,\UU)$ is a unitary asymptote of $\TT$.
\begin{itemize}
\item[\textup{(a)}] The linear mapping $\Gamma= \Gamma_{(X,\UU)}\co \{\UU\}'\to \I(\TT,\UU), Y'\mapsto Y'X$ is invertible.
\item[\textup{(b)}] The unitary asymptote $(X,\UU)$ has optimal norm-control with $\k_{\textup{op}}= \|\Gamma^{-1}\|$.
\end{itemize}
\end{pro}

\noindent
\textbf{Proof.}
The mapping $\Gamma$ is bounded linear, for any unitary intertwining pair of $\TT$.
It is obvious that $\Gamma$ is bijective, whenever $(X,\UU)$ is a unitary asymptote of $\TT$.
In this case $\Gamma^{-1}$ is also bounded by the Open Mapping Theorem.

Let $(X',\UU')$ be a unitary intertwining pair of $\TT, \; X'=Y'X$ with $Y'\in \I(\UU,\UU')$.
Since $|Y'|\in \{\UU\}'$, it follows that $\|Y'\|= \| |Y'| \| \le \|\Gamma^{-1}\|\cdot \| |Y'| X\|$.
In view of Proposition \ref{generating} we know that $\| |Y'| X\| =\| X'\|$, and so $\|Y'\|\le \|\Gamma^{-1}\|\cdot \|X'\|$.
It is clear that $\|\Gamma^{-1}\|$ is the optimal norm-control.

\rightline{$\square$}

The question of the uniqueness of unitary asymptotes can be settled easily.

\begin{pro}
\label{unique}
Suppose that $(X,\UU)$ is a unitary asymptote of $\TT$, and $(\wt X,\wt\UU)$ is a unitary intertwining pair of $\TT$.
\begin{itemize}
\item[\textup{(a)}] Then $(\wt X,\wt\UU)$ is a unitary asymptote of $\TT$ if and only if $(\wt X,\wt\UU) \approx (X,\UU)$.
\item[\textup{(b)}] If $Z\in\I(\wt\UU,\UU)$ is invertible with $Z\wt X=X$, then 
\[\k_{\rm op}(\wt X,\wt\UU) \le \|ÄZ\|\cdot \k_{\rm op}(X,\UU) \le \|Z\|\cdot \|Z^{-1}\|\cdot  \k_{\rm op}(\wt X,\wt\UU).\]
\end{itemize}
\end{pro}

\noindent
\textbf{Proof.}
Let us assume that $(X,\UU)$ and $(\wt X,\wt\UU)$ are unitary asymptotes of $\TT$.
Then there exist unique transformations $Z\in \I(\wt\UU,\UU)$ and $\wt Z\in \I(\UU,\wt\UU)$ satisfying the conditions $Z\wt X=X$ and $\wt Z X=\wt X$.
Since $I\cdot X= X= Z\wt X= Z\wt Z X$ and $I, Z\wt Z\in \{\UU\}'$, it follows that $I=Z\wt Z$.
The equality $I=\wt Z Z$ can be obtained similarly.
Thus, $Z$ is invertile with $\wt Z = Z^{-1}$, and so $(X,\UU)\approx (\wt X, \wt\UU)$.

If $(\wt X,\wt\UU) \approx (X,\UU)$, then there exists an invertible $Z\in\I(\wt\UU,\UU)$ with $Z\wt X=X$.
Given any unitary intertwining pair $(X',\UU')$ of $\TT$, there is a unique $Y'\in\I(\UU,\UU')$ such that $X'=Y'X$. 
Then $\wt Y'= Y'Z\in \I(\wt\UU,\UU')$ is the unique transformation satisfying $\wt Y'\wt X= X'$. 
Thus, $(\wt X, \wt\UU)$ is a unitary asymptote of $\TT$.
Furthermore, if $(X,\UU)$ has norm-control with $\k$, then $\|\wt Y'\|\le \|Y'\|\cdot \|Z\|\le \k \|Z\| \|X'\|$, and so $(\wt X,\wt\UU)$ has norm-control with $\k\|Z\|$.
Therefore, $\k_{\rm op}(\wt X,\wt\UU) \le \|Z\|\cdot \k_{\rm op}(X,\UU)$.
Similarly, starting with $\wt X= Z^{-1}X$ we obtain $\k_{\rm op}(X,\UU)\le \|Z^{-1}\|\cdot \k_{\rm op}(\wt X,\wt\UU)$.

\rightline{$\square$}

\begin{rem}
\label{absolute-control}
\textup{
Let $(X,\UU)$ be a unitary asymptote of $\TT$.
Clearly, $\k_{\rm op}(X,\UU)=0$ holds exactly, when $X=0$.
Let us assume that $X\ne 0$.
Given any $0\ne c\in\C$, and applying Proposition \ref{unique} with $\wt\UU=\UU$ and $Z=cI$, we obtain that $(\wt X=(1/c)X,\UU)$ is also a unitary asymptote of $\TT$ and $\k_{\rm op}(\wt X,\UU)= |c| \k_{\rm op}(X,\UU)$.
Hence the products $\|\wt X\|\cdot \k_{\rm op}(\wt X,\UU)= \|X\|\cdot \k_{\rm op}(X,\UU)$ are independent of the particular choice of $c$; we call this common value the \emph{absolute optimal norm-control} of $(X,\UU)$, and denote it by $\k_{\rm aop}(X,\UU)$.
Clearly, $\k_{\rm op}(X,\UU) = \|X\|^{-1} \k_{\rm aop}(X,\UU)$.
Notice also that $\|\Gamma_{(X,\UU)}\|\le \|X\|$ yields that $\k_{\rm op}(X,\UU)= \|\Gamma_{(X,\UU)}^{-1}\|\ge 1/\|X\|$, whence $\k_{\rm aop}(X,\UU)\ge 1$ follows.
}
\end{rem}

\begin{pro}
\label{unique2}
If $(X,\UU)$ and $(\wt X,\wt\UU)$ are unitary asymptotes of $\TT$ with norm-control $1$, and $\|X\|\le 1, \|\wt X\|\le1$, then $(X,\UU)\eq (\wt X,\wt \UU)$; furthermore $\|X\|=1$, provided $X\ne 0$.
\end{pro}

\noindent
\textbf{Proof.}
Let us consider the previous transformations $Z, \wt Z$.
Since $\|Z\|\le \|X\|\le1$ and $\|Z^{-1}\|= \|\wt Z\|\le \|\wt X\|\le 1$, it follows that $Z$ is unitary.
Thus $\|Z\|=1$, provided $Z$ does not act between the zero spaces.

\rightline{$\square$}

Similarity preserves the existence of unitary asymptotes.

\begin{pro}
\label{similarity}
If $(\wt X, \wt\UU)$ is a unitary asymptote of the commuting $n$-tuple $\wt\TT$ with norm-control $\k$ and $ Z\in\I(\TT, \wt\TT)$ is invertible, then $(\wt XZ,\wt\UU)$ is a unitary asymptote of $\TT$ with norm-control $\k\|Z^{-1}\|$.
Furthermore, we have $\k_{\rm aop}(\wt X Z,\wt\UU)\le \|Z\|\cdot\|Z^{-1}\|\cdot \k_{\rm aop}(\wt X,\wt\UU)$.
\end{pro}

\noindent
\textbf{Proof.}
Let $(X',\UU')$ be a unitary intertwining pair of $\TT$.
Then $(X'Z^{-1},\UU')$ is a unitary intertwining pair of $\wt\TT$, and so there is a unique $Y'\in\I(\wt\UU,\UU')$ such that $X'Z^{-1}= Y'\wt X$.
The last equality is equivalent to $X'= Y'(\wt X Z)$.
Furthermore, $\|Y'\|\le \k \|X'Z^{-1}\|\le \k \|Z^{-1}\| \|X'\|$.

\rightline{$\square$}

As an application, we are able to provide simple examples for unitary asymptotes.
 
\begin{cor}
\label{simple-ua}
Let $\WW$ be a commuting $n$-tuple of unitaries on a non-zero space, and let $Z\in\I(\TT,\WW)$ be invertible.
Then
\begin{itemize}
\item[\textup{(a)}] $\; (I,\WW)$ is a unitary asymptote of $\WW$ with optimal norm-control $1$;
\item[\textup{(b)}] $\; (Z,\WW)$ is a unitary asymptote of $\TT$ with norm-control $\|Z^{-1}\|$.
\end{itemize}
\end{cor}

The following example shows that the norm-control with $\|Z^{-1}\|$ is not optimal (even not comparable to the optimal), in general.
For simplicity, we deal with the single operator case; i.e., $n=1$ is assumed.

\begin{ex}
\label{2dimensional}
\textup{
Given any $\l\in\T\setminus\{1\}$, let us consider the operators
\[
T= \left[
\begin{matrix}
\l & 0 \\
1 & 1
\end{matrix}
\right] \quad \hbox{ and } \quad 
W= \left[
\begin{matrix}
\l & 0 \\
0 & 1
\end{matrix}
\right], \]
acting on the Hilbert space $\C^2$.
It can be easily checked that
\[
Z=\left[
\begin{matrix}
1 & 0 \\
\frac{1}{1-\l} & 1
\end{matrix}
\right] \in \I(T,W) \; \hbox{ is invertible with } \; Z^{-1}=
\left[
\begin{matrix}
1 & 0 \\
\frac{-1}{1-\l} & 1
\end{matrix}
\right].\]
Thus, $(Z,W)$ is a unitary asymptote of $T$ with norm-control $\|Z^{-1}\|$.
It is clear that $\|Z^{-1}\|\ge 1/|1-\l|$ can be made as large as we wish, approaching $\l$ to $1$.
On the other hand, in view of Proposition \ref{control}, the optimal norm-control for $(Z,W)$ is $\k_{\textup{op}}= \|\Gamma^{-1}\|$, where $ \Gamma\co \{W\}' \to \I(T,W), Y' \mapsto Y'Z$.
The general form of $Y'$ is 
\[
Y'=
\left[
\begin{matrix}
\eta_1 & 0 \\
0 & \eta_2
\end{matrix}
\right], \quad \hbox{ whence } \; Y'Z =
\left[
\begin{matrix}
\eta_1 & 0 \\
\frac{\eta_2}{1-\l} & \eta_2
\end{matrix}
\right].
\]
Obviously, 
\begin{eqnarray*}
\|\Gamma^{-1}\| & = & \sup\left\{ \frac{\|Y'\|}{\|Y'Z\|}: 0\ne Y'\in\{W\}'\right\} \\
& = & \left(\inf\{\|Y'Z\|: Y'\in \{W\}', \|Y'\|=1\}\right)^{-1}.
\end{eqnarray*}
If $1= \|Y'\| = \max\{|\eta_1|, |\eta_2|\}$, then $\|Y'Z\| \ge \max\{|\eta_1|, |\eta_2|\}=1$; moreover, $\|Y'Z\|=1$ when $\eta_1=1$ and $\eta_2=0$.
Thus, $\k_{\textup{op}}= \|\Gamma^{-1}\|=1$.
}
\end{ex}

Necessary conditions for the existence of non-zero unitary intertwining pairs can be given in terms of the spectral radius.
We say that the commuting $n$-tuple $\TT=(T_1,\dots,T_n)$ is invertible, if each $T_i$ is invertible ($i\in\N_n$); then $\TT^{-1}:= (T_1^{-1},\dots,T_n^{-1})$ is also a commuting $n$-tuple.

\begin{pro}
\label{sp-radius}
Let us assume that there is a unitary intertwining pair $(X,\UU)$ of $\TT$, with $X\ne0$.
\begin{itemize}
\item[\textup{(a)}] Then $r(\TT^\kk)\ge 1$ holds, for all $\kk\in\Z_+^n$.
\item[\textup{(b)}] If $\TT$ is invertible, then $r(\TT^\kk)=1$ is true, for all $\kk\in\Z^n$.
\end{itemize}
\end{pro}

\noindent
\textbf{Proof.}
For any $\kk\in \Z_+^n$ and $j\in \N$, we have $X\TT^{j\kk}=\UU^{j\kk}X$, whence it follows that $(\TT^{j\kk})^*X^*X\TT^{j\kk} = X^*X$.
Thus $\|X^*X\|\le \|\TT^{j\kk}\|^2\cdot \|X^*X\|$, where $\|X^*X\|= \|X\|^2>0$.
We infer that $\|\TT^{j\kk}\|\ge 1$, what implies that $r(\TT^\kk)= \lim_{j\to\i} \|(\TT^\kk)^j\|\ge1$.

If $\TT$ is invertible, then $X\TT^{j\kk}= \UU^{j\kk}X$ holds for every $\kk\in\Z$ and $j\in\N$.
As before, we obtain $r(\TT^\kk)\ge1$ and also $r(\TT^{-\kk})\ge1$; hence $r(\TT^\kk)=1$.

\rightline{$\square$}

We say that $\TT$ is of \emph{$0$-type}, if $X=0$ whenever $(X,\UU)$ is a unitary intertwining pair of $\TT$.
In this case $(0,\textbf{0})$ is a degenerate unitary asymptote of $\TT$, where $\textbf{0}=(0,\dots,0)$ acts on the zero space $\{0\}$.
As a consequence, we obtain further examples for unitary asymptotes.

\begin{cor}
\label{zero-type}
If $\TT$ is invertible and $r(\TT^\kk)\ne 1$ for some $\kk\in\Z^n$, then \textup{ $(0, \textbf{0})$} is a unitary asymptote of $\TT$.
\end{cor}

We turn to orthogonal sums.
If $\TT_j=(T_{j1},\dots,T_{jn})$ is a commuting $n$-tuple of operators on $\HH_j$, for $j=1,2$, then $\TT_1\op \TT_2:= (T_{11}\op T_{21},\dots,T_{1n}\op T_{2n})$ is a commuting $n$-tuple of operators on the Hilbert space $\HH_1\op\HH_2$.

\begin{pro}
\label{orthogonal}
For $j=1,2$, let $\TT_j=(T_{j1},\dots,T_{jn})$ be a commuting $n$-tuple of operators on the Hilbert space $\HH_j$, and let us assume that $(X_j,\UU_j)$ is a unitary asymptote of $\TT_j$ with norm-control $\k_j$.
Then $(X=X_1\op X_2, \UU= \UU_1\op \UU_2)$ will be a unitary asymptote of $\TT= \TT_1\op \TT_2$ with norm-control $\k= \sqrt{2} \max(\k_1,\k_2)$.
\end{pro}

\noindent
\textbf{Proof.}
It is clear that $(X,\UU)$ is a minimal unitary intertwining pair of $\TT$.
Suppose that $(X',\UU')$ is also a unitary intertwining pair of $\TT$.
For $j=1,2$, the restriction $X_j':= X'|\HH_j$ belongs to $\I(\TT_j,\UU')$, and so there is a unique $Y_j'\in \I(\UU_j,\UU')$ such that $X_j'= Y_j'X_j$.
Then $Y'= [Y_1'\;\; Y_2']$ is the unique transformation in $\I(\UU,\UU')$ satisfying $Y'X= X'$.
Furthermore, we have $\|Y'\|\le \sqrt{2} \max(\|Y_1'\|, \|Y_2'\|) \le \sqrt{2} \max(\k_1\|X_1'\|, \k_2\|X_¡'\|) \le \sqrt{2} \max(\k_1,\k_2) \|X'\|.$

\rightline{$\square$}

The next statement deals with the opposite direction.

\begin{pro}
\label{orthogonal2}
Let us assume that $(X,\UU)$ is a unitary asymptote of $\TT=\TT_1\op\TT_2$, acting on $\HH=\HH_1\op \HH_2$, with norm-control $\k$.
For $j=1,2$, let us consider the reducing subspace $\K_j= \vee \{\UU^\kk X\HH_j : \kk\in\Z^n\}$ of $\UU$, and the restrictions $\UU_j:= \UU|\K_j$ and $X_j:= X|\HH_j\in\I(\TT_j,\UU_j)$.
Then $(X_j,\UU_j)$ will be a unitary asymptote of $\TT_j\; (j=1,2)$ with norm-control $\k$; furthermore, $(X,\UU) \approx (X_1\op X_2, \UU_1\op\UU_2)$.
\end{pro}

\noindent
\textbf{Proof.}
It is clear that $(X_j,\UU_j)$ is a minimal unitary intertwining pair of $\TT_j\; (j=1,2)$.
Suppose that $(X_j',\UU_j')$ is also a unitary intertwining pair of $\TT_j$.
Then $(\wt X_j,\UU_j')$ will be a unitary intertwining pair of $\TT$ with $\wt X_j h:= X_j'P_jh\; (h\in\HH)$, where $P_j\in\L(\HH)$ denotes the orthogonal projection onto $\HH_j$.
Hence there exists unique $\wt Y_j\in\I(\UU,\UU_j')$ such that $\wt Y_jX=\wt X_j$.
We infer that $Y_j':= \wt Y_j|\K_j\in \I(\UU_j,\UU_j')$ and $Y_j'X_j= \wt YX|\HH_j = \wt X_j|\HH_j = X_j'$.
Assuming that $Y_j''\in\I(\UU_j,\UU_j')$ and $Y_j''X_j=X_j'$, the relations $Y_j''\UU_j^\kk X_j= {\UU_j'}^\kk Y_j'' X_j= {\UU_j'}^\kk Y_j'X_j = Y_j' \UU_j^\kk X_j \; (\kk\in \Z^n)$ yield that $Y_j''= Y_j'$.
Finally, $\|Y_j'\|\le \|\wt Y_j\| \le \k \|\wt X_j\|= \k \|X_j'\|$, and so $(X_j, \UU_j)$ is a unitary asymptote of $\TT_j$ with norm-control $\k$.

The last statement is an immediate consequence of Propositions \ref{orthogonal} and \ref{unique}.

\rightline{$\square$}

Concluding this section we provide example for an operator, which does not have a unitary asymptote (and so it is necessarily not of $0$-type).
Furthermore, we shall see that, taking restrictions to invariant subspaces or taking adjoints, the existence of unitary asymptotes can be lost.

\begin{ex}
\label{no-ua}
\textup{
Let $\{e_j\}_{j\in\Z}$ be an orthonormal basis in the Hilbert space $\wt\E$, and let $\wt S\in \L(\wt\E)$ be the bilateral shift defined by $\wt Se_j=e_{j+1}\; (j\in\Z)$.
Let us consider the operator $\wt R= 2\wt S$ and its restriction $R= \wt R|\E$ to the invariant subspace $\E= \vee\{e_j\}_{j\in\Z_+}$.
It is evident that $(\wt X,\wt S)$ is a unitary intertwining pair of $R$, where the injective transformation $\wt X$ is defined by $\wt X e_j= 2^{-j} e_j\; (j\in \Z_+)$.
Thus, $R$ is far from being of $0$-type.}

\textup{Let us assume that there exists a unitary asymptote $(X,U)$ of $R$.
For any $\l\in\T, \; (X_\l,\l I_\C)$ is a unitary intertwining pair of $R$, where $X_\l h = \sum_{j=0}^\i 2^{-j}\l^j \la h, e_j\ra\; (h\in \E)$.
Hence there is a unique $Y_\l \in\I(U,\l I_\C)$ such that $X_\l =Y_\l X$.
We may easily infer that $\l$ is an eigenvalue of $U$.
Since the eigenspaces corresponding to distinct eigenvalues are orthogonal to each other, it follows that $\sum_{\l\in\T}\op \ker(U-\l I)$ is a subspace of the domain $\K$ of $U$, which means that the Hilbert space $\K$ is not separable.
But this is impossible, because $\E$ is separable and $(X,U)$ is minimal (see Proposition \ref{minimal}).
Therefore,} $R$ does not have a unitary asymptote.

\textup{
The operator $\wt R$ is invertible and $r(\wt R)=2$, hence $\wt R$ is of $0$-type by Corollary \ref{zero-type}.
Let $W$ be a unitary operator on a non-zero Hilbert space $\F$.
Then $(0\op I, 0\op W)$} is a unitary asymptote of 
\textup{$\wt T= \wt R\op W$ (see Proposition \ref{orthogonal}).
The subspace $\HH= \E\op\F$ is invariant for $\wt T$; but the} restriction $T:= \wt T|\HH= R\op W$ does not have a unitary asymptote,
\textup{ since $R$ fails to have a unitary asymptote (see Proposition  \ref{orthogonal2}).}

\textup{
If $(X,U)$ is a unitary intertwining pair of $R^*$, then $0= \|X(R^*)^{j+1} e_j\| = \| U^{j+1}Xe_j\| = \|Xe_j\| \; (j\in\Z_+)$; hence $X=0$, and so $R^*$ is of $0$-type.
We conclude that the} adjoint $T^*= R^*\op W^*$ has a unitary asymptote: $(0\op I, 0\op W^*)$.
\end{ex}

Analogous counterexamples can be obtained for commuting $n$-tuples with $n>1$, observing that $(X,U)$ is a unitary asymptote of $T$ exactly when the pair $(X, (U,I,\dots,I))$  is a unitary asymptote of $(T,I,\dots,I)$.

\section{Orbit conditions}
\label{orbits}

Let $\TT=(T_1,\dots,T_n)$ be a commuting $n$-tuple of operators on the Hilbert space $\HH$.
The \emph{orbit-infimum} of $\TT$ at the vector $h\in\HH$ is defined by 
$\oinf(\TT,h):= \inf\{\|\TT^\kk h\|: \kk\in\Z_+^n\}$.
If $(X',\UU')$ is any unitary intertwining pair of $\TT$, then  $\|X'h\|= \|{\UU'}^\kk X'h\|= \|X'\TT^\kk h\| \le \|X'\|\cdot \|\TT^\kk h\|\; (h\in\HH, \kk\in \Z_+^n)$, and so the \emph{Upper Orbit Condition} (UOC) automatically holds:
\[\|X' h\|\le \|X'\| \cdot \oinf(\TT,h)\quad \hbox{ for all }\; h\in\HH.\]
The opposite \emph{Lower Orbit Condition} (LOC) \emph{ with bound-control} $\k\in\R_+$ is an extra requirement:
\[\oinf(\TT,h)\le \k \|X' h\| \quad \hbox{ for all } \; h\in\HH.\]
The latter condition together with minimality are sufficient for being a unitary asymptote.

\begin{pro}
\label{orbit-condition}
Let $(X,\UU)$ be a minimal unitary intertwining pair of $\TT$.
If \textup{ (LOC) } holds for $(X,\UU)$ with bound-control $\k$, then $(X,\UU)$ is a unitary asymptote of $\TT$ with norm-control $\k$.
\end{pro}

\noindent
\textbf{Proof.}
Let $(X',\UU')$ be any unitary intertwining pair of $\TT$.
Given any $h\in\HH$, the (UOC) for $(X',\UU')$ yields $\|X'h\|\le \|X'\|\cdot \oinf(\TT,h)$, while the (LOC) for $(X,\UU)$ results $\oinf(\TT,h)\le \k \|Xh\|$.
Thus $\|X'h\|\le \|X'\| \cdot \k \|Xh\|\; (h\in\HH)$, and so there is a unique $Y_+'\in \L((X\HH)\ba, (X'\HH)\ba)$ such that $Y_+'Xh= X'h \; (h\in \HH)$;
in particular, $\|Y_+'\|\le \k \|X'\|$.
By the minimality of $(X,\UU)$, there exists a unique $Y'\in \L(\K, \K')$ such that $Y'\UU^\kk y= {\UU'}^\kk Y_+' y$ holds for all $y\in (X\HH)\ba$ and $\kk\in\Z^n$.
It can be easily seen that $Y'\in \I(\UU,\UU'), Y'X=X'$ and $\|Y'\|= \|Y_+'\|\le \k \|X'\|.$

\rightline{$\square$}

Let us consider the homogeneous set $\HH_{00}(\TT):= \{h\in\HH: \oinf(\TT,h)=0\}$.
If $(X',\UU')$ is a unitary intertwining pair of $\TT$, then the (UOC) implies that $\HH_{00}(\TT)\subset \ker X'$.
Assuming that $(X,\UU)$ is a unitary asymptote of $\TT$, the relation $X'=Y'X$ yields that $\ker X \subset \ker X'$.
In particular, if $(\wt X,\wt\UU)$ is another unitary asymptote of $\TT$, then $\ker\wt X= \ker X$; this common nullspace $\HH_0(\TT)$ is called the \emph{annihilating subspace} of $\TT$.
Clearly, $\HH_0(\TT)$ contains the set $\HH_{00}(\TT)$.
It is immediate that coincidence is ensured by the (LOC).

\begin{pro}
\label{coincidence}
If the unitary asymptote $(X,\UU)$ of $\TT$ satisfies the \textup{ (LOC)}, then $\HH_0(\TT)= \HH_{00}(\TT)$.
\end{pro}

Let us assume that $\TT$ has a unitary asymptote.
Relying on the annihilating subspace $\HH_0(\TT)$, we may classify the $n$-tuples:
\begin{itemize}
\item $\TT\in C_{0 \cdot}$, if $\HH_0(\TT)=\HH$; this happens exactly, when $(0, \mathbf{0})$ is a unitary asymptote, that is when $\TT$ is \emph{of $0$-type};
\item $\TT\in C_{*\cdot}$, if $\HH_0(\TT)\ne \HH$, then $\TT$ is called \emph{asymptotically non-vanishing};
\item $\TT\in C_{1\cdot}$, if $\HH_0(\TT)= \{0\}\ne\HH$; then $\TT$ is called \emph{asymptotically strongly non-vanishing}.
\end{itemize}
If the adjoint $\TT^*$ also has a unitary asymptote, then $\TT\in C_{\cdot j}$ by definition if $\TT^*\in C_{j\cdot}\; (j=0,*,1)$.
Finally, $C_{ij}= C_{i\cdot}\cap C_{\cdot j}$.

Clearly, $\HH_{00}(\TT)=\HH$ implies that $\TT\in C_{0\cdot}$.
However, the following example shows that the $n$-tuples in the class $C_{0\cdot}$ can be separated form $0$.
It shows also that the  (LOC) \emph{ is not necessary for the existence of unitary asymptotes.}

\begin{ex}
\label{no-LOC}
\textup{Let us assume that $\TT$ is a commuting $n$-tuple of invertible operators on the non-zero Hilbert space $\HH$, and that $\|T_i^{-1}\|<1$ for every $1\le i\le n$.
By Corollary \ref{zero-type}, $\TT$ is of $0$-type, that is $\HH_0(\TT)=\HH$ and $\TT\in C_{0\cdot}$.
On the other hand, for every $0\ne h\in\HH$ we have $\oinf(\TT,h)= \|h\|>0$, even more $\lim_\kk\|\TT^\kk h\|=\i$.
Hence $\HH_{00}(\TT)=\{0\}$ and the (LOC) fails.
}
\end{ex}

The (LOC) does hold in the power bounded setting.
Let us assume that $\TT$ is a commuting $n$-tuple of power bounded operators on $\HH$, that is we have $s_i= \sup\{\|T_i^{k_i}\|: k_i\in\Z_+\}<\i$  for every $1\le i\le n$.
Then $\|\TT^\kk\|\le s_1\cdots s_n$ is true, for all $\kk\in\Z_+^n$,
and so $s:= \sup\{\|\TT^\kk\|: \kk\in\Z_+\}<\infty$.

Let $M\co \ell^\i(\Z_+^n)\to \C$ be an \emph{invariant mean}, that is $M$ is a linear functional with $\|M\|=1=M(\mathbf{1})$, and $M(\xi_\mathbf{t})= M(\xi)$ holds for every $\xi\in\ell^\i(\Z_+)$ and $\mathbf{t}\in\Z_+^n$, where $\xi_{\mathbf{t}}(\kk):= \xi(\kk+\mathbf{t})\; (\kk\in\Z_+^n)$.
We shall use the notation: $\Mlim_\kk \xi(\kk)= M(\xi)$.
(See, e.g., \cite{Day}  and \cite{Patterson}.)
The bounded, sesquilinear functional $w(x,y):= \Mlim_\kk\la \TT^\kk x,\TT^\kk y\ra \; (x,y\in\HH)$ can be represented by a unique positive operator $A\in\L(\HH)$ in the sense that $w(x,y)= \la Ax,y\ra$ holds, for all $x,y\in\HH$.
Let us consider the Hilbert space $\K_+= (A\HH)\ba$, and the bounded, linear transformation $X_+\co \HH\to \K_+,\; h\mapsto A^{1/2}h,$  with $\|X_+\|^2=\| A\|= \|w\|\le s^2$.
For any $i\in\N_n$, for every $h\in\HH$, we have $\|X_+h\|^2= \la Ah,h\ra = \Mlim_\kk \la \TT^\kk h, \TT^\kk h\ra = \Mlim_\kk \la \TT^\kk T_i h, \TT^\kk T_i h\ra = \la AT_i h, T_ih\ra = \|X_+ T_i h\|^2$.
Hence, there exists a unique isometry $V_i\in \L(\K_+)$ such that $V_iX_+= X_+T_i$.
The relations $V_iV_jX_+= V_iX_+T_j = X_+T_iT_j = X_+T_jT_i= V_jV_iX_+$ yield that $V_iV_j= V_jV_i$.
The commuting $n$-tuple of isometries $\mathbf{V}= (V_1,\dots,V_n)$ can be extended to a commuting $n$-tuple of unitaries $\UU=(U_1,\dots,U_n)$ acting on a larger space $\K$ such that $\mathbf{V}= \UU|\K_+$, and $\K= \vee\{\UU^\kk \K_+: \kk\in\Z^n\}$ (see Proposition I.6.2 in [NFBK]).
We readily obtain that $(X,\UU)$ is a \emph{minimal unitary intertwining pair} of $\TT$, where $X\co \HH\to\K, h\mapsto X_+h$, and $\|X\|= \|X_+\|\le s$.
Since
\[\|Xh\|^2 = \|X_+ h\|^2= \la Ah,h\ra = \Mlim_\kk \|\TT^\kk h\|^2\ge (\oinf(\TT,h))^2\]
holds, for every $h\in\HH$, it follows that $(X,\UU)$ \emph{satisfies the} (LOC) \emph{with bound-control $1$.}

In view of Propositions \ref{orbit-condition}, \ref{coincidence} and \ref{unique}, and applying also the inequality $\|\TT^{\kk + \mathbf{l}}h\|\le s \|\TT^\kk h\| \; (\kk, \mathbf{l}\in \Z_+^n, h\in\HH)$, we obtain the following proposition.

\begin{pro}
\label{power-bounded}
Let $\TT= (T_1,\dots,T_n)$ be a commuting $n$-tuple of power bounded operators on $\HH$, with $s= \sup\{\|\TT^\kk\|: \kk\in\Z_+^n\}<\i$.
For any invariant mean $M\co \ell^\i(\Z_+^n)\to\C$, let $(X_M,\UU_M)$ denote the pair, constructed in the previous canonical way.
\begin{itemize}
\item[\textup{(a)}] Then $(X_M,\UU_M)$ is a unitary asymptote of $\TT$ with $\|X_M\|\le s$, satisfying \textup{(LOC)} with bound-control $1$, and so $\k_{\rm aop}(X_M,\UU_M)\le s$.
\item[\textup{(b)}] If $(X_{\wt M}, \UU_{\wt M})$ is the unitary asymptote of $\TT$ corresponding to another invariant mean $\wt M$, then there is a unique invertible $Z\in \I(\UU_M, \UU_{\wt M})$ such that $X_{\wt M} = ZX_M, \|Z\|\le s$ and $\|Z^{-1}\|\le s$
\item[\textup{(c)}] The annihilating subspace of $\TT$ consists of the stable vectors of $\TT$:
\[\HH_0(\TT)= \HH_{00}(\TT)= \left\{h\in\HH: \lim_\kk \|\TT^\kk h\|=0\right\}.\]
\end{itemize}
\end{pro}

We recall that every operator with bounded orbits is power bounded by the Uniform Boundedness Principle.
It would be interesting to give example for a commuting $n$-tuple $\TT$ of operators, such that $\TT$  is not power bounded, and $\TT$ has a unitary asymptote $(X,\UU)$ satisfying (LOC).
An example like that would provide negative answer to the following question.

\begin{que}
\label{necessary-condition}
Is power boundedness of $\TT$ a necessary condition for the existence of a unitary asymptote $(X,\UU)$ (possibly with injective $X$), satisfying \textup{(LOC)}?
\end{que}

\begin{rem}
\label{contractions}
\textup{ Let us assume now that $\TT= (T_1,\dots,T_n)$ is a commuting $n$-tuple of} contractions,
\textup{that is $\|T_i\|\le 1$ for each $i\in\N_n$.
Then, for every $h\in\HH, \{\|\TT^\kk h\|\}_{\kk\in \Z_+^n}$ is a decreasing generalized sequence (net) with respect to the natural inductive partial ordering on $\Z_+^n$, and so it is convergent: $\lim_\kk\|\TT^\kk h\|= \oinf(\TT,h)$.
The polar identity shows that $\{\la\TT^\kk x, \TT^\kk y\ra\}_\kk$ also converges, for every $x,y\in\HH$.
Hence all invariant means $M$ result in the same unitary asymptote $(X_M, \UU_M)= (X, \UU)$, with 
$\|Xh\|= \lim_\kk \|\TT^\kk h\|\; (h\in \HH)$.
}
\end{rem}

\begin{rem}
\label{components}
\textup{
(a) Let us assume again that $\TT= (T_1,\dots,T_n)$ is a commuting $n$-tuple of power bounded operators on $\HH$.
For any $i\in\N_n,\; \oinf(\TT,h)\le \oinf(T_i,h)$ holds for every $h\in \HH$, and so $\HH_0(T_i)= \HH_{00}(T_i)\subset \HH_{00}(\TT)= \HH_0(\TT)$.
In particular, if $\TT$ is asymptotically (strongly) non-vanishing, then $T_i$ is also asymptotically (strongly) non-vanishing.}

\textup{ (b) The following example shows that the converse implication is not true. 
Let us consider tha analytic Toeplitz operator $T_\f f= \f f$ on the Hardy--Hilbert space $H^2$, where $\f\in H^\i, \|\f\|=1$ and $\O(\f)= \{\z\in\T: |\f(\z)|=1\}$ and $\T\setminus\O(\f)$ are of positive arcwise measure.
Then $U_\f g= \f g$ is a unitary operator on the Hilbert space $L^2(\O(\f))$, and $X_\f\co H^2\to L^2(\O(\f)), f \mapsto f|\O(\f)$ is a contractive transformation with dense range.
It can be easily verified that $(X_\f, U_\f)$ is a minimal unitary intertwining pair of $T_\f$, and that the (LOC) holds with bound-control $1$.
Thus, $(X_\f,U_\f)$ is a unitary asymptote of $T_\f$ with norm-control $1$.
Furthermore, the F. \&  M. Riesz Theorem implies that $\HH_0(T_\f)=0$, that is $T_\f\in C_{1\cdot}$.
(In connection with the theory of Hardy spaces we refer to \cite{Hof}.)
}

\textup{ Let us choose $\f_1, \f_2\in H^\i$ so that $\|\f_1\|_\i = \|\f_2\|_\i=1$, and $\O(\f_1), \O(\f_2)$ are disjoint sets of positive measure.
Let us consider the commuting pair of contractions $\TT= (T_{\f_1}, T_{\f_2})$ on $H^2$.
We have seen that $T_{\f_1}$ and $T_{\f_2}$ are asymptotically strongly non-vanishing.
On the other hand, $\lim_{j\to\i} \|T_{\f_1}^j T_{\f_2}^j h\|=0$ for every $h\in \HH$, and so $\HH_0(\TT)= H^2$, that is $\TT\in C_{0\cdot}$.}
\end{rem}

The following proposition, inspired by a remark of Maria Gamal', shows that unitary asymptotes, not satisfying (LOC), arise naturally.

\begin{pro}
\label{inverse}
Let $\TT$ be a commuting $n$-tuple of invertible operators on $\HH$, and let $(X,\UU)$ be a unitary asymptote of $\TT$ with norm-control $\k$.
\begin{itemize}
\item[\textup{(a)}] Then $(X,\UU^{-1})$ is a unitary asymptote of $\TT^{-1}$, with norm-control $\k$.
\item[\textup{(b)}] If $\TT$ is power bounded, then the following conditions are equivalent:
\begin{itemize}
\item[\textup{(i)}] $\TT^{-1}$ has a unitary asymptote satisfying \textup{(LOC)};
\item[\textup{(ii)}] the given $(X,\UU^{-1})$ satisfies \textup{(LOC)};
\item[\textup{(iii)}] $X$ is invertible;
\item[\textup{(iv)}] $\TT$ is similar to a commuting $n$-tuple of unitaries.
\end{itemize}
\end{itemize}
\end{pro}

\noindent
\textbf{Proof.}
(a): Let $(X',\UU'^{-1})$ be a unitary intertwining pair of $\TT^{-1}$.
Since $(X',\UU')$ is a unitary intertwining pair of $\TT$, there is a unique $Y'\in \I(\UU,\UU')$ such that $X'=Y'X$.
Then $Y'\in\I(\UU^{-1},\UU'^{-1})$ also holds, whence the statement follows.
We note also that $(X\HH)\ba$ is reducing for $\UU$, and so $X$ has dense range by minimality.

(b): Suppose that $\TT$ is power bounded: $s:= \sup\{\|\TT^\kk\|: \kk\in \Z_+^n\}<\infty$.
(Then the existence of the unitary asymptote is guaranteed by Proposition \ref{power-bounded}.)
If $\TT^{-1}$ has a unitary asymptote satisfying (LOC), then $(X,\UU^{-1})$ satisfies (LOC) too by similarity (see Proposition \ref{unique}).

For every $h\in\HH$, we have $\oinf(\TT^{-1},h)\ge \|h\|/s$.
Indeed, if $\oinf(\TT^{-1},h)< \|h\|/s$, then $\|\TT^{-\kk} h\|<\|h\|/s$ holds for some $\kk\in\Z_+^n$, and so we infer that $\|h\|= \|\TT^\kk \TT^{-\kk}h\|\le s \|\TT^{-\kk}h\|< \|h\|$, what is a contradiction.
Therefore, if (LOC) holds for $(X,\UU^{-1})$ with a bound-control $\k>0$,
then the inequalities
\[\frac{1}{s} \|h\| \le \oinf(\TT^{-1},h) \le \k \|Xh\| \quad (h\in\HH)\]
imply that $X$ is invertible, and so $\TT$ is similar to the unitary $n$-tuple $\UU$.

Finally, if $\TT$ is similar to a unitary $n$-tuple, then $\TT^{-1}$ is power bounded, and so it has a unitary asymptote satisfying (LOC) by Proposition \ref{power-bounded}.

\rightline{$\square$}

Now we give a specific example for an invertible $\TT\in C_{1\cdot}$, where (LOC) fails for $\TT^{-1}$.

\begin{ex}
\label{LOCfails}
\textup{Let us consider the Hilbert space $\ell^2(\Z^n,\b)$ consisting of the sequences $\xi\co \Z^n\to\C$ satisfying $\|\xi\|_\b^2:= \sum_{\kk\in\Z^n} |\xi(\kk)|^2 \b(\kk)^2<\infty$, where $\log_2\b(\kk)= - \sum_{i=1}^n \min(0,k_i)$ for $\kk=(k_1,\dots,k_n)\in \Z^n$.
The commuting invertible contractions $\TT=(T_1,\dots,T_n)$ are defined on $\ell^2(\Z^n,\b)$ by $(T_i \xi)(\kk)= \xi(\kk-\mathbf{e}^{(i)})$, where $\mathbf{e}^{(i)}= (\d_{i1},\dots,\d_{in})\; (i\in\N_n)$.
The commuting unitaries $\UU=(U_1,\dots,U_n)$ are defined on $\ell^2(\Z^n)$ by the same formula: $(U_i\xi)(\kk)= \xi(\kk-\mathbf{e}^{(i)})$.
The transformation
$X \co \ell^2(\Z^n,\b)\to \ell^2(\Z^n), \; \xi\mapsto \xi$
intertwines $\TT$ with $\UU$, and 
$\|X\xi\|= \lim_\kk \|\TT^\kk\xi\|_\b = \|\xi\|\le \|\xi\|_\b$ holds, for every $\xi\in\ell^2(\Z^n,\b)$.
We obtain that $(X,\UU)$ is a unitary asymptote of $\TT$, satisfying (LOC) with bound-control 1.
Since the injective $X$ is not invertible, we infer by Proposition \ref{inverse} that
 (LOC) fails for the unitary asymptote $(X,\UU^{-1})$ of $\TT^{-1}\in C_{1\cdot}$.}
\end{ex}

\section{Quasianalytic $n$-tuples of operators}
\label{quasianalytic}

Let $\TT= (T_1,\dots,T_n)$ be a commuting $n$-tuple of operators on the Hilbert space $\HH$, and let $\UU= (U_1,\dots,U_n)$ be a commuting $n$-tuple of unitaries on $\K$.
\emph{Let us assume that $(X,\UU)$ is a unitary asymptote of $\TT$ with norm-control $\k$, where $X\ne 0$, and so $\TT\in C_{*\cdot}$.}

Given any $C\in \{\TT\}', \; (XC,\UU)$ is a unitary intertwining pair of $\TT$.
Hence, there exists a unique $D\in\I(\UU,\UU)= \{\UU\}'$ such that $XC=DX$.
The transformation
\[\g=\g_X\co \{\TT\}' \to \{\UU\}', C\mapsto D\]
is called the \emph{commutant mapping} associated with the unitary asymptote $(X,\UU)$.
Obviously, $\g(\TT^\kk)= \UU^\kk$ holds for all $\kk\in\Z_+^n$.

\begin{pro}
\label{commutant-mapping}
The commutant mapping $\g$  is a bounded algebra-homomor-phism, with $\|\g\|\le \k \|X\|$; taking $\k$ optimal, we obtain $\|\g\|\le \k_{\rm aop}(X,\UU)$.
\end{pro}

\noindent
\textbf{Proof.}
Given any $C_1,C_2\in \{\TT\}'$, we have $C_1C_2\in \{\TT\}'$ and $\g(C_1C_2)X= X(C_1C_2)= (XC_1)C_2= \g(C_1)XC_2= \g(C_1)\g(C_2)X$.
Since $\g(C_1C_2), \g(C_1)\g(C_2)$ belong to $ \{\UU\}'$, it follows that $\g(C_1C_2)= \g(C_1)\g(C_2)$.
The linearity of $\g$ can be verified similarly.
Finally, because of the norm-control, we infer $\|\g(C)\|\le \k \|XC\|\le \k \|X\|\cdot\|C\|$ for all $C\in \{\TT\}'$.

\rightline{$\square$}

\begin{rem}
\label{cm-related}
\textup{ Assuming that $(\wt X, \wt\UU)$ is another unitary asymptote of $\TT$, let $Z\in\I(\UU,\wt\UU)$ be the invertible transformation satisfying $ZX=\wt X$ (see Proposition \ref{unique}).
Given any $C\in\{\TT\}'$, let us consider the operators $D=\g_X(C)\in \{\UU\}'$ and $\wt D= \g_{\wt X}(C)\in\{\wt\UU\}'$.
Since $(X,\UU)$ is a unitary asymptote of $\TT,\; \wt DZX= \wt D \wt X= \wt X C= ZXC= ZDX\in \I(\TT,\wt\UU)$ and $\wt D Z, ZD \in \I(\UU,\wt\UU)$, it follows that $\wt DZ=ZD$.
Observing that $(Z,\wt\UU)$ is a unitary asymptote of $\UU$ (see Corollary \ref{simple-ua}), we can see that $\wt D= ZDZ^{-1}= \g_Z(D)$.
Therefore, we obtain that $\g_{\wt X}= \g_Z \circ\g_X$.
}
\end{rem}

The \emph{hyperinvariant subspace lattice} $\Hlat \TT$ of $\TT$ consists of those subspaces $\M$ of $\HH$, which are invariant for the commutant $\{\TT\}'$, that is $C\M\subset\M$ holds for all $C\in\{\TT\}'$.
The relation $\{\TT\}'\subset \{T_i\}'$ readily implies that $\Hlat\TT \supset \Hlat T_i$ is true, for every $i\in\N_n$.
For illustration we give an example.

\begin{ex}
\label{Hlat-ex}
\textup{Let $S$ be the unilateral shift on the Hardy--Hilbert space $H^2$; that is $Sf=\h f\; (f\in H^2)$, where $\h(z)=z$ is the identical function.
Let us consider the commuting pair $\TT= (T_1,T_2)$ of the operators $T_1= S\op S$ and $T_2= I\op 0$ acting on the Hilbert space $\HH= H^2\op H^2$.}

\textup{It is well-known that the commutant of $S$ consists of the analytic Toeplitz operators on $H^2$ (see, e.g., Section 147 in \cite{Halmos}).
Hence the invariant subspaces are also hyperinvariant for $S$, and by Beurling's theorem they are of the form $\th H^2$, where $\th\in H^\i$ is an inner function ($|\th(\z)|=1$ for a.e. $\z\in\T$) or $\th=0$.
Since $C=[C_{ij}]_2\in \{T_1\}'$ exactly when $C_{ij}\in \{S\}'$ for all $1\le i,j\le 2$, it can be easily seen that $\Hlat T_1= \{\th(H^2\op H^2): \th\in H^\i$ inner or $\th=0\}$.
Obviously, $\{T_2\}' = \{A\op B: A, B\in \L(H^2)\}$ and $\Hlat T_2 = \{\{0\}, H^2\op \{0\}, \{0\}\op H^2, \HH\}$.
We conclude that $\{\TT\}' = \{T_1\}'\cap \{T_2\}' = \{T_\f \op T_\psi: \f, \psi\in H^\i\}$, and so $\Hlat \TT= \{\th H^2 \op \eta H^2: \th, \eta \in H^\i$ inner or $0\}$.}
\end{ex}

Since $\{\UU\}'$ is a $C^*$-algebra, the subspaces in $\Hlat \UU$ are reducing.
Furthermore, spectral theory provides many subspaces in $\Hlat \UU$, as we shall see soon.
The hyperinvariant subspaces of $\UU$ induce hyperinvariant subspaces of $\TT$.

\begin{pro}
\label{induced-HSP}
For every \textup{ $\NN\in\Hlat \UU$}, we have \textup{ $\M= X^{-1}\NN\in \Hlat\TT$}; in particular, \textup{$\HH_0(\TT)= X^{-1}(\{0\})\in \Hlat \TT$}.
\end{pro}

\noindent
\textbf{Proof.}
Given any $C\in\{\TT\}'$, let us consider $D=\g(C)\in \{\UU\}'$.
For any $h\in\M, \; Xh\in\NN$ implies $DXh\in\NN$.
Since $XCh= DXh$, it follows that $Ch\in\M$.

\rightline{$\square$}

Now we turn to the spectral analysis of the commuting $n$-tuple $\UU$ of unitaries, recalling some basic facts.
Let $C^*[\UU]= \vee\{\UU^\kk: \kk\in\Z^n\}$ be the abelian $C^*$-algebra, generated by $\UU$ in $\L(\K)$.
(We remind the reader that $\K\ne \{0\}$.)
Let $\SS(\UU)$ denote the spectrum of $C^*[\UU]$, consisting of all non-zero complex homomorphisms $\Lambda\co C^*[\UU]\to \C$.
The spectrum $\SS(\UU)$ is a compact Hausdorff space with the topology, induced by the weak-$*$ topology.
The Gelfand transform $\Gamma_\UU\co A \mapsto \widehat A$, defined by $\widehat A(\Lambda):= \Lambda(A)$, is a $C^*$-algebra isomorphism from $C^*[\UU]$ onto the $C^*$-algebra $C(\SS(\UU))$ of all continuous functions on $\SS(\UU)$.
Since $C^*[\UU]$ is generated by $\UU$, it follows that $\wt\tau_\UU\co \Lambda\mapsto \Lambda(\UU):= (\Lambda(U_1),\dots, \Lambda(U_n))$ is an injective, continuous mapping from $\SS(\UU)$ into the $n$-dimensional torus $\T^n$.
The range of $\wt\tau_\UU$, denoted by $\s(\UU)$, is called the (joint) \emph{spectrum} of $\UU$; it is a non-empty, compact subset of $\T^n$.
Since the spectrum with respect to a $C^*$-subalgebra is the same as in $\L(\K)$, it can be easily seen that $\s(U_i)= \{\Lambda(U_i)= \pi_i(\Lambda(\UU)): \Lambda\in \SS(\UU)\}= \pi_i(\s(\UU))$, where $\pi_i\co \C^n\to \C, (z_1,\dots,z_n)\mapsto z_i$ is the projection on the $i$-th coordinate ($i\in\N_n$).
We note that $\s(\UU)$ coincides with the (joint) \emph{approximate point spectrum} $\s_{\textup{ap}}(\UU)$, consisting of all $ (\l_1,\dots,\l_n)\in \C^n$ satisfying the condition $ \lim_{j\to\i} \sum_{i=1}^n \|U_i x_j- \l_i x_j\|= 0$ with a sequence of unit vectors $\{x_j\}_{j\in\N}$ in $\K$ (see \cite{Lyubich} and \cite{Muller}).

The homeomorphism $\tau_\UU\co \SS(\UU) \to \s(\UU),\; \Lambda \mapsto \Lambda(\UU)$ induces the isomorphism $\tau_\UU^\#\co C(\s(\UU)) \to C(\SS(\UU)), f \mapsto f \circ \tau_\UU$.
The composition $\PP_\UU= \Gamma_\UU^{-1} \circ \tau_\UU^\#$ is a faithful representation of the $C^*$-algebra $C(\s(\UU))$ on $\K$.
By the Spectral Theorem (see, e.g., \cite{ConwayFA}) there exists a unique spectral measure $E_{\UU,0}\co \B_{\s(\UU)}\to \P(\K)$ such that $\PP_\UU(f)= \int_{\s(\UU)} f \, dE_{\UU,0}$, for all $f\in C(\s(\UU))$.
(Here $\B_{\s(\UU)}$ denotes the system of Borel subsets of $\s(\UU)$, and $\P(\K)$ stands for the set of orthogonal projections in $\K$.)
The mapping $E_{\UU,0}$ can be extended to the spectral measure $E_\UU$, defined on the $\s$-algebra $\B_n$ of all Borel subsets of $\T^n$ by $E_\UU(\o):= E_{\UU,0}(\o\cap\s(\UU))$, and called as the \emph{asymptotic spectral measure} of $\TT$, corresponding to $(X,\UU)$.
The support of $E_\UU$ is the spectrum $\s(\UU)$; i.e., $\T^n\setminus\s(\UU)$ is the largest open subset of $\T^n$ with $E_\UU(\T^n\setminus\s(\UU))=0$.
Clearly, $U_i= \PP_\UU(\pi_i|\s(\UU))= \int_{\T^n} z_i \, dE_\UU(\zz)$ is true for every $i\in\N_n$, and so $\UU^\kk= \int_{\T^n} \zz^\kk\, dE_\UU(\zz)$ holds for all $\kk\in\Z^n$.
For any $i\in\N_n$, the spectral measure of $U_i$ can be derived from $E_\UU$ by $E_{U_i}(\o)= E_\UU(\pi_i^{-1}(\o)\cap\T^n)$, where $\o\in \B_1= \B_\T$.
Conversely, if $\o= \o_1 \times \cdots \times \o_n$ with $\o_1,\dots,\o_n\in \B_1$, then $E_\UU(\o)= E_{U_1}(\o_1)\cdots E_{U_n}(\o_n)$.
We note yet that $D\in \{\UU\}'$ if and only if $D$ commutes with every spectral projection $E_\UU(\o) \; (\o\in \B_n)$.
Hence, \emph{the spectral subspaces $E_\UU(\o)\K$ are all hyperinvariant for $\UU$.}

The next step towards quasianalicity is to take localization.
The \emph{localization} of the spectral measure $E=E_\UU$ at the vector $v\in\K$ is the finite positive Borel measure $E_v$ on $\T^n$, defined by $E_v(\o):= \la E(\o)v,v\ra \; (\o\in\B_n).$
Clearly, $E_v(\T^n)= \|v\|^2$.
Let $M_+(\T^n)$ be the set of all finite, positive Borel measures on $\T^n$.
Given any $\m, \nu\in M_+(\T^n)$, the notation $\m\ac\nu$ means that $\m$ is absolutely continuous with respect to $\nu$, that is $\nu(\o)=0$ implies $\m(\o)=0\; (\o\in\B_n)$.
The measures $\m$ and $\nu$ are equivalent, in notation: $\m\ae\nu$, if $\m\ac\nu$ and $\nu\ac\m$, that is $\m(\o)=0$ holds exactly when $\nu(\o)=0$.

Let $\wt\UU= (\wt U_1,\dots, \wt U_n)$ be also a commuting $n$-tuple of unitary operators, acting on the non-zero Hilbert space $\wt\K$, and let $\wt E\co  \B_n\to\P(\wt\K)$ be the spectral measure of $\wt\UU$.
Given any $Z\in\I(\UU,\wt\UU)$, let us consider the unitary intertwining pair $(\wt X, \wt\UU)$ of $\TT$, where $\wt X:= ZX$.
The following technical lemma relates localizations of $E$ and $\wt E$.

\begin{lem}
\label{related-loc}
If $Z$ is injective, then $\wt E_{\wt X h}= E_{|Z| Xh} \ae E_{Xh}$ hold, for every $h\in \HH$.
\end{lem}

\noindent
\textbf{Proof.}
Taking the polar decomposition $Z= W |Z|$, it can be easily seen that $|Z|\in \{\UU\}'$ is injective and $W\in\I(\UU,\wt\UU)$ is an isometry.
The subspace $\wt\K_1=W\K$ is reducing for $\wt\UU$, so the restriction $\wt\UU_1:= \wt\UU|\K_1$ is a commuting $n$-tuple of unitaries on $\wt\K_1$.
Let $\wt E_1\co \B_n\to\P(\wt\K_1)$ be the spectral measure of $\wt\UU_1$.
Clearly, $W_1\in \I(\UU, \wt\UU_1)$ is a unitary transformation, where $W_1 v= Wv \; (v\in\K)$.
Relying on the uniqueness part of the Riesz Representation Theorem, it is easy to show that $W_1E(\o)W_1^*= \wt E_1(\o)= \wt E(\o)|\wt\K_1$ hold for all $\o\in\B_n$.
Thus, for any $h\in \HH$, we have
\begin{eqnarray*}
\wt E_{\wt X h}(\o) &=& \la\wt E(\o)\wt X h, \wt X h\ra = \la \wt E_1(\o)W_1 |Z| Xh, W_1 |Z| Xh\ra \\
& = & \la W_1 E(\o) |Z| Xh, W_1 |Z| Xh\ra = \la E(\o) |Z| Xh, |Z| Xh\ra \\
&=& E_{|Z|Xh}(\o) \quad (\o\in\B_n).
\end{eqnarray*}

Since $|Z|\in \{\UU\}'$, it follows that $|Z| E(\o)= E(\o) |Z|$, and so $E(\o) |Z| Xh = |Z| E(\o)Xh$.
Taking into account that $|Z|$ is injective, we infer that $E_{|Z| Xh}(\o)=0$ if and only if $E_{Xh}(\o)=0$.
(Notice that $E_v(\o)= \|E(\o)v\|^2$, for all $v\in\K$.)
Therefore, we obtain that $\wt E_{\wt X h}= E_{|Z| Xh} \ae E_{Xh}$.

\rightline{$\square$}

We say that the \emph{unitary intertwining pair} $(X',\UU')$ of $\TT$ is \emph{quasianalytic}, if $X'\ne0$ and $E'_{X'u} \ae E'_{X'v}$ holds for all $u,v\in \HH\setminus \{0\}$, where $E'\co \B_n \to \P(\K')$ denotes the spectral measure of $\UU'$.
The \emph{commuting $n$-tuple} $\TT$ of operators is called \emph{quasianalytic}, if it has a quasianalytic unitary asymptote $(X,\UU)$.
In view of Proposition \ref{unique} and Lemma \ref{related-loc}, we can see that if one of the unitary asymptotes is quasianalytic, then so are all of them.

\begin{pro}
\label{asymp-behav}
If $\TT$ is a quasianalytic commuting $n$-tuple of operators, then it is asymptotically strongly non-vanishing: $\TT\in C_{1\cdot}$.
\end{pro}

\noindent
\textbf{Proof.}
Let $(X,\UU)$ be a unitary asymptote of $\TT$, and let $u\in\HH$ be a vector such that $Xu\ne0$.
Then, for every $0\ne v\in\HH$, the relation $E_{Xv}\ae E_{Xu}\ne0$ yields that $E_{Xv}\ne0$, and so $Xv\ne0$.
Thus, $\HH_0(\TT)= \ker X = \{0\}$.

\rightline{$\square$}

Quasianalycity is a property of homogeneity type.
If it is broken, then hyperinvariant subspaces arise.

\begin{thm}
\label{proper-HSP}
If the asymptotically non-vanishing, commuting $n$-tuple $\TT$ of operators is not quasianalytic, then \textup{$\Hlat \TT$} is non-trivial.
\end{thm}

\noindent
\textbf{Proof.}
Since $\TT\in C_{*\cdot}$, it has a unitary asymptote $(X,\UU)$, where $\ker X =\HH_0(\TT) \ne \HH$.
In view of Proposition \ref{induced-HSP}, we may assume that $\HH_0(\TT)= \{0\}$, when $X$ is injective.
Let $E\co \B_n\to \P(\K)$ be the spectral measure of $\UU$.
Since $\TT$ is not quasianalytic, we can find non-zero vectors $u,v\in\HH$ so that $E_{Xu} \not\ac E_{Xv}$.
Thus, there exists $\o_1\in\B_n$ such that $E_{Xv}(\o_1)=0$ and $E_{Xu}(\o_1)>0$.
Setting $\o_2= \T^n\setminus \o_1$, let us consider the decomposition $\K= \NN_1\op \NN_2$, where $\NN_j= E(\o_j)\K\in \Hlat\UU$ for $j=1,2$.
We know from Proposition \ref{induced-HSP} that $\M= X^{-1}\NN_2\in \Hlat\TT$.
The equality $E_{Xv}(\o_1)=0$ yields that $Xv\in\NN_2$, whence $v\in\M$ follows.
On the other hand, $E_{Xu}(\o_1)>0$ implies that $Xu\not\in\NN_2$, and so $u\not\in\M$.
Consequently, $\M$ is a proper hyperinvariant subspace of $\TT$.

\rightline{$\square$}

\begin{cor}
\label{proper-HSP2}
Let $\TT$ be a commuting $n$-tuple of class $C_{*\cdot}$.
If there exists an injection $Y\in \I(\wt\UU,\TT)$, where $\wt\UU$ is a commuting $n$-tuple of unitaries on $\wt\K$ and the spectrum $\s(\wt\UU)$ is not a single point of $\T^n$, then $\TT$ is not quasianalytic, and so \textup{$\Hlat\TT$} is non-trivial.
\end{cor}

\noindent
\textbf{Proof.}
Let $\wt E$ denote the spectral measure of $\wt\UU$.
Since $\s(\wt\UU)$ is not a singleton, a Borel set $\o_1\in\B_n$ can be given so that $\wt E(\o_1)\ne 0$ and $\wt E(\o_2)\ne 0$ with $\o_2= \T^n\setminus\o_1$.
Choosing non-zero vectors $v_j\in \wt E(\o_j)\wt\K\; (j=1,2)$, the relation $\wt E(\o_1) \wt E(\o_2)= 0$ yields that the non-zero localizations $\wt E_{v_1}$ and $\wt E_{v_2}$ are singular to each other.

Let $(X,\UU)$ be a unitary asymptote of $\TT$, and let $E$ be the spectral measure of $\UU$.
Suppose that $\TT$ is quasianalytic.
Then $X$ is injective by Proposition \ref{asymp-behav}, and so $Z=XY\in \I(\wt\UU,\UU)$ is also injective.
Since $(I,\wt\UU)$ is a unitary asymptote of $\wt\UU$ (see Corollary \ref{simple-ua}), an application of Lemma \ref{related-loc} results in that $E_{Zv_j} \ae \wt E_{v_j}$ for $j=1,2$.
Thus $E_{X(Yv_1)}$ is not equivalent to $E_{X(Yv_2)}$, which means that $\TT$ is not quasianalytic.
We arrived at a contradiction, and so $\TT$ cannot be quasianalytic.
Theorem \ref{proper-HSP} implies that $\Hlat \TT$ is non-trivial.

\rightline{$\square$}

We conclude this section with an example.

\begin{ex}
\label{qa-example}
\textup{ For every $i\in\N_n$, let $\f_i\in H^\i$ be given so that $\|\f_i\|_\i=1$ and $\O(\f_i)= \{\z\in\T: |\f_i(\z)|=1\}$ is of positive Lebesgue measure.
Let us assume also that $\O= \cap_{i=1}^n \O(\f_i)$ and $\T\setminus\O$ are of positive measure.
Let us consider the commuting $n$-tuple $\TT= (T_{\f_1},\dots, T_{\f_n})$ of analytic Toeplitz operators on the Hardy--Hilbert space $H^2$; that is $T_{\f_i}f= \f_i f$ for all $f\in H^2\; (i\in\N_n)$.
Let us consider also the commuting $n$-tuple $\UU=( U_{\f_1},\dots,U_{\f_n})$ of unitary multiplication operators on the Hilbert space $L^2(\O)= \h_\O L^2(\T)$, where $\h_\O$ stands for the characteristic function of $\O$ and $L^2(\T)$ is defined with respect to the normalized Lebesgue measure $m$ on $\T$.
Thus, $U_{\f_i}g= \f_i g$ for all $g\in L^2(\O)\; (i\in\N_n)$.
It is clear that $X\in\I(\TT,\UU)$, where $Xf:= \h_\O f\; (f\in H^2)$.
In view of the F. \& M. Riesz Theorem $ f(\z)\ne0$ for a.e. $\z\in\T$, whenever $0\ne f\in H^2$.
Hence $X$ is injective.
Furthermore, for every $0\ne f\in H^2, \; Xf$ is cyclic for the unitary operator $U_\O$, defined by $U_\O g= \h g \; (g\in L^2(\O))$; that is $\vee_{j\in\Z_+} U_\O^j Xf= L^2(\O)$ (see, e.g., Section 146 in \cite{Halmos} and \cite{Bram}).
In particular, $X$ has dense range, and so $(X,\UU)$ is a minimal unitary intertwining pair of $\TT$.
Since the (LOC) holds with bound-conrol $\k=1$, it follows that $(X,\UU)$ is a unitary asymptote of $\TT$ with norm-control $\k=1$ (see Proposition \ref{orbit-condition}).
We obtain also that $\HH_0(\TT) = \ker X = \{0\}$, that is $\TT\in C_{1\cdot}$.
Taking into account that $U_\O\in \{\UU\}'$, we conclude that $Xf$ is cyclic for the commutant $\{\UU\}': \; \vee\{DXf: D\in\{\UU\}'\}= L^2(\O)$, for every $0\ne f\in H^2$.
Therefore, $\TT$ must be quasianalytic; see the proof of Theorem \ref{proper-HSP}}.
\end{ex}

\section{Local and global spectral invariants}
\label{sp-invariants}

Let $\TT= (T_1,\dots,T_n)$ be again a commuting $n$-tuple of operators on the separable Hilbert space $\HH$.
Let $\UU= (U_1,\dots,U_n)$ be a commuting $n$-tuple of unitaries on the Hilbert space $\K$, and let us assume that $(X,\UU)$ is a unitary asymptote of $\TT$  with $X\ne0$; that is $\TT\in C_{*\cdot}$.
Let $E\co \B_n\to \P(\K)$ denote the spectral measure of $\UU$, where $\B_n$ stands for the $\s$-algebra of all Borel subsets of $\T^n$.

The Borel sets $\o_1, \o_2\in\B_n$ \emph{essentially coincide with respect to} $E$, in notation: $\o_1 \eqE \o_2$, if $E(\o_1\triangle \o_2)=0$, that is when $E(\o_1)= E(\o_2)$.
This extended identity is clearly an equivalence relation on $\B_n$.
The Borel set $\o_1$ is \emph{essentially contained in $\o_2$ with respect to} $E$, in notation: $\o_1\subE \o_2$, if $E(\o_1\setminus\o_2)=0$, that is when $E(\o_1)\le E(\o_2)$.
This relation is a partial ordering on $\B_n$ with the extended identity (otherwise speaking, on the  quotient set $\B_n/E$).

\begin{pro}
\label{clattice-E}
The pair $(\B_n,\subE)$ is a complete lattice.
\end{pro}

\noindent
\textbf{Proof.}
Given any $\o_1, \o_2\in\B_n$, let $\O:= \o_1\cup\o_2$.
It is clear that $\o_1, \o_2\subE\O$.
Suppose that $\o_1, \o_2 \subE \o'\in\B_n$.
Since $\O\setminus \o'\subset (\o_1\setminus\o')\cup (\o_2\setminus\o')$, it follows that $E(\O\setminus\o')\le E((\o_1\setminus \o')\cup (\o_2\setminus\o'))= E(\o_1\setminus\o')+ E(\o_2\setminus\o')- E(\o_1\setminus\o')E(\o_2\setminus\o')=0.$
Hence $\O\subE \o'$, and so $\O= \o_1 \veeE \o_2$.
(We could have written $\O\in \o_1 \veeE\o_2$, because $\o_1\veeE \o_2$ is determined only up to sets of $E$-measure $0$.)

Given any sequence $\{\o_l\}_{l\in\N}$ in $\B_n$, we obtain by induction that $\O_N:= \cup_{l=1}^N \o_l= \veeE\{\o_l\}_{l=1}^N$ for every $N\in\N$.
Let us consider the Borel set $\ol\o:= \cup_{l\in\N}\o_l= \cup_{N\in\N}\O_N$.
Clearly, $\ol\o \overset{E}\supset \o_l$ holds for all $l\in\N$.
Suppose that $\B_n\ni \o'\overset{E}\supset \o_l$ for all $l\in\N$.
Then, for every $N\in\N$, $\o'\overset{E}\supset \O_N$, that is $E(\o')\ge E(\O_N)$.
Since $E(\O_N)$ converges to $E(\ol\o)$ in the strong operator topology (as $N\to\i$), it follows that $E(\o')\ge E(\ol\o)$, that is $\o' \overset{E}\supset \ol\o$.
Therefore, $\ol\o= \veeE\{\o_l\}_{l\in\N}$.

Let us consider now an arbitrary system $\{\o_j\}_{j\in J}$ of Borel sets on $\T^n$.
Since the Hilbert space $\HH$ is separable, it follows by the minimality of the unitary asymptote that $\K$ is also separable.
Hence the subspace $\NN= \vee_{j\in J} E(\o_j)\K$ is separable too.
Let $\{y_k\}_{k\in\N}$ be a dense sequence in $\NN$.
For every $k\in\N$, there exists a countable subset $J_k$ of $J$ such that $y_k\in \vee_{j\in J_k} E(\o_j)\K$.
Then the set $\wt J:= \cup_{k\in\N} J_k$ is also countable: $\wt J= \{j_l\}_{l\in\N}$.
The relations $\NN\supset \vee_{l\in\N} E(\o_{j_l})\K \supset \vee\{y_k\}_{k\in\N}=\NN$ imply that $\NN= \vee_{l\in\N} E(\o_{j_l})\K= E(\ol\o)\K$, where $\ol\o:= \cup_{l\in\N}\o_{j_l}\in\B_n$.
For every $j\in J, \; E(\ol\o)\K=\NN \supset E(\o_j)\K$ yields $\ol\o  \overset{E}\supset \o_j$.
Suppose that $\o' \overset{E}\supset \o_j$ for all $j\in J$.
Then $E(\o')\K\supset E(\o_j)\K$ for all $j\in J$, whence $E(\o')\K\supset \NN= E(\ol\o)\K$ follows, and so $\o' \overset{E}\supset \ol\o$.
Consequently, $\ol\o=\veeE \{\o_j\}_{j\in J}$.

For every Borel sets $\o, \o'\in \B_n, \; \o \subE\o'$ is equivalent to $ (\o')^c \subE \o^c$, valid for the comlements.
Thus, $\overset{E}\wedge\{\o_j: j\in J\}= \left( \veeE\{\o_j^c: j\in J\}\right)^c$ also exists.

\rightline{$\square$}

Let us consider the localization of $E$ at a vector $y\in\K: E_y(\o)= \la E(\o)y,y\ra = \|E(\o)y\|^2 \; (\o\in\B_n)$.
It is immediate that $E_y(\o)=0$ if and only if $E(\o)y=0$, and $E_y(\o)= \|y\|^2$ exactly when $E(\o)y=y$.
Furthermore, $E_y(\o)= \|y\|^2$ if and only if $E_y(\o^c)=0$.
The Borel set
\[\o(\UU,y):= \overset{E}\wedge\{\o\in\B_n: E_y(\o)= \|y\|^2\} = \left( \veeE\{\o'\in\B_n: E_y(\o')=0\}\right)^c\]
is called the \emph{local residual set} of $\UU$ at $y$; it is determined up to the extended identity with respect to $E$.
Clearly,
\[E(\o(\UU,y))\K = \cap\{E(\o)\K: y\in E(\o)\K\}\]
and
\[E(\o(\UU,y)^c)\K= \vee\{E(\o')\K: E(\o')y=0\}.\]
The next proposition contains further basic properties; the proof is left to the reader.

\begin{pro}
\label{basic-pr}
Using the previous notation, we have:
\begin{itemize}
\item[\textup{(a)}] For every $y\in\K, \; E(\o(\UU,y))\K$ is the smallest spectral subspace containing $y$.
\item[\textup{(b)}] For any $y_1, y_2\in\K$, we have $E_{y_1}\ae E_{y_2}$ if and only if $\o(\UU,y_1)= \o(\UU,y_2)$.
\item[\textup{(c)}] $\veeE \{\o(\UU,y): y\in\K\} = \T^n$.
\end{itemize}
\end{pro}

For every vector $h\in\HH, \; \o(\TT,h):= \o(\UU,Xh)$ is called the \emph{local residual set} of $\TT$ at $h$.
This Borel set is determined up to $E$-identity, and it is independent of the particular choice of the unitary asymptote (see Lemma \ref{related-loc} and its proof).
Since the reducing subspace $\ol\K= \vee \{E(\o(\TT,h))\K: h\in\HH\}$ contains the range of $X$, we infer that $\ol\K=\K$, and so
\[\veeE\{\o(\TT,h): h\in\HH\}= \T^n.\]
The Borel set $\pi(\TT):= \overset{E}\wedge\{\o(\TT,h): h\in \HH\setminus\{0\}\}$ is called the \emph{quasianalytic spectral set} of $\TT$.
It is also determined up to $E$-identity and independent of the special choice of $(X,\UU)$.
In view of Proposition \ref{basic-pr} we obtain the following theorem.

\begin{thm}
\label{qa-coincidence}
The commuting $n$-tuple $\TT\in C_{*\cdot}$ is quasianalytic if and only if $\pi(\TT)= \T^n$.
\end{thm}

Quasianalycity can be also characterized in terms of cyclicity.
To deduce this theorem we need some auxiliary results concerning the spectral analysis of $\UU$.
Though these must be known, we sketch the proofs, since we could not find direct references.

The measure $\mu=\m_\UU\in M_+(\T^n)$ is called a \emph{scalar spectral measure} of $\UU$, if it is equivalent to the spectral measure $E: \m(\o)=0$ if and only if $E(\o)=0\; (\o\in\B_n)$.

\begin{pro}
\label{scalar-spm}
We have:
\begin{itemize}
\item[\textup{(a)}] For any $y\in\K$, the localization $E_y$ is a scalar spectral measure of $\UU$ if and only if $\o(\UU,y)=\T^n$.
\item[\textup{(b)}] There exists a vector $y\in\K$ such that $\o(\UU,y)=\T^n$.
\end{itemize}
\end{pro}

\noindent
\textbf{Proof.}
(a): For any $\o\in\B_n, E_y(\o)=0$ if and only if $\o \subE \o(\UU,y)^c$, and $E(\o)=0$ exactly when $\o \eqE \emptyset$.
Hence $E_y$ is equivalent to $E$ if and only if $\o(\UU,y)^c\eqE\emptyset$, that is when $\o(\UU,y) \eqE \T^n$.

(b): It is easy to construct a sequence $\{y_l\}_{l\in\N}$ of vectors in $\K$, such that $\|y_l\|<2^{-l}$ for all $l\in\N$, the spectral projections $\{E(\o_l)\}_{l\in\N}$ are orthogonal to each other, and $\vee_{l\in\N}E(\o_l)\K=\K$, where $\o_l= \o(\UU,y_l)$.
Let us consider the vector $y:= \sum_{l=1}^\i y_l$.
Since the spectral projections commute, we infer that the systems $\{y_l\}_{l\in\N}$ and $\{E(\o(\UU,y))y_l\}_{l\in\N}$ are orthogonal.
Hence, the equalities $y= E(\o(\UU,y))y= \sum_{l=1}^\i E(\o(\UU,y))y_l$ imply that $ \sum_{l=1}^\i \|y_l\|^2 = \|y\|^2 = \sum_{l=1}^\i \|E(\o(\UU,y))y_l\|^2$, where $\|y_l\|^2\ge \|E(\o(\UU,y))y_l\|^2\; (l\in\N)$.
Thus, for every $l\in\N, \; E(\o(\UU,y))y_l=y_l$, whence $\o(\UU,y) \overset{E}\supset \o_l$ follows.
We conclude that $\o(\UU,y)\overset{E}\supset \veeE \{\o_l\}_{l\in\N}=\T^n$.

\rightline{$\square$}

Let $\B(\T^n)$ denote the $C^*$-algebra of all bounded, Borel measurable functions on $\T^n$.
The extended functional calculus
\[\FF_{\UU,0}\co \B(\T^n) \to \L(\K),\; g\mapsto g(\UU):= \int_{\T^n} g \, dE\]
for $\UU$ is a $*$-homomorphism, that is a representation of $\B(\T^n)$ on $\K$.
Let $\m\in M_+(\T^n)$ be a scalar spectral measure of $\UU$.
Since for every $y\in\K, \; \|g(\UU)y\|^2= \la |g|^2(\UU)y,y\ra = \int_{\T^n}|g|^2\, dE_y$ where $E_y\ac\m$, and since $E_y\ae \m$ holds for some $y\in\K$ by Proposition \ref{scalar-spm}, it follows that $g(\UU)=0$ if and only if $g(\zz)=0$  for almost every $\zz\in\T^n$ with respect to $\m$.
Thus, the well-defined mapping
\[\FF_\UU\co L^\i(\m)\to \L(\K),\; g\mapsto g(\UU)= \int_{\T^n}g\, dE\]
is an injective $*$-homomorphism, that is a faithful representation of the $C^*$-algebra $L^\i(\m)$, and so it is isometric.
Identifying $L^\i(\m)$ with the dual of the Banach space $L^1(\m)$, we may consider the weak-$*$ topology on $L^\i(\m)$.
Since, for every $y\in\K, \; \la g(\UU)y,y\ra = \int_{\T^n}g\, dE_y= \int_{\T^n} g \f_y\, d\m$ where $0\le \f_y= dE_y/d\m\in L^1(\m)$, we obtain that $\FF_\UU$ is continuous also with respect to the weak-$*$ topology on $L^\i(\m)$ and the weak operator topology on $\L(\K)$.

It is clear that $\{\UU\}''= C^*[\UU]''$ holds for the bicommutants, and by the Double Commutant Theorem this algebra coincides with the closure of $C^*[\UU]$ in the weak operator topology (see Section 12  in \cite{ConwayOT}).
Thus, $\{\UU\}''$ is the abelian von Neumann algebra, generated by $\UU$.

\begin{pro}
\label{sp-subspaces}
The range of the functional calculus $\FF_\UU$ coincides with $\{\UU\}''$, and so \textup{$\Hlat\UU= \{E(\o)\K: \o\in\B_n\}$}.
\end{pro}

\noindent
\textbf{Proof.}
The closed unit ball $(L^\i(\m))_1= \{g\in L^\i(\m): \|g\|_\i\le1\}$ is weak-$*$ compact by the Banach--Alaoglu Theorem.
Thus its image $\FF_\UU(L^\i(\m))_1= (\FF_\UU L^\i(\m))_1$ is compact in the weak operator topology, and so it is weak-$*$ closed in $\L(\K)$.
Applying the Krein--Smulyan Theorem we obtain that the $C^*$-algebra $\FF_\UU L^\i(\m)$ is weak-$*$ closed, hence it is also closed in the weak operator topology (see Section V.12 in \cite{ConwayFA} and Section I.3.4 in \cite{Dixmier}).
Therefore  $\FF_\UU L^\i(\m)$ is an abelian von Neumann algebra, containing $C^*[\UU]= \PP_\UU C(\T^n)= \FF_\UU C(\T^n)$.
Since $C(\T^n)$ is weak-$*$ dense in $L^\i(\m)$ by Lusin's theorem, we conclude that $\FF_\UU L^\i(\m)$ is the closure of $C^*[\UU]$ in the weak operator toõpology.
Consequently, $\FF_\UU L^\i(\m)= \{\UU\}''$.

We have already mentioned that the spectral subspaces $E(\o)\K\; (\o\in\B_n)$ are all hyperinvariant for $\UU$.
Let us assume now that $\NN\in\Hlat \UU$.
Then the orthogonal projection $P_\NN$ onto $\NN$ belongs to the bicommutant $\{\UU\}''$.
Hence, there exists a function $g\in L^\i(\m)$ such that $P_\NN= g(\UU)$.
Since $g^2(\UU)= g(\UU)^2= g(\UU)$, it follows that $g^2=g$, and so $g= \h_\o$ holds with some $\o\in \B_n$.
We infer that $\NN= P_\NN \K= \h_\o(\UU)\K= E(\o)\K$.

\rightline{$\square$}

\begin{thm}
\label{ciclicity}
For the commuting $n$-tuple $\TT\in C_{*\cdot}$, the following conditions are equivalent:
\begin{itemize}
\item[\textup{(i)}] $\TT$ is quasianalytic;
\item[\textup{(ii)}] the vector $Xh$ is cyclic for the commutant $\{\UU\}'$, for every non-zero $h\in\HH$;
\item[\textup{(iii)}] for every non-zero $h\in\HH, \; E(\o)Xh\ne0$ whenever $E(\o)\ne0 \; (\o\in\B_n)$.
\end{itemize}
\end{thm}

\noindent
\textbf{Proof.}
Suppose that, for a non-zero $h\in\HH$, the vector $Xh$ is not cyclic for $\{\UU\}'$.
Then the hyperinvariant subspace $\NN= \{DXh: D\in \{\UU\}'\}\ba$ of $\UU$, induced by $Xh$, is not equal to $\K$.
By Proposition \ref{sp-subspaces}, there exists a Borel set $\o_h\in\B_n$ such that $E(\o_h)\K=\NN$.
Since $\o(\TT,h)= \o_h$ does not essentially coincide with $\T^n$ with respect to $E$, we infer that $\TT$ is not quasianalytic (see Proposition \ref{scalar-spm} and Theorem \ref{qa-coincidence}).
Thus, (i) implies (ii).

Since the spectral subspaces are hyperinvariant, it is immediate that (ii) yields (iii).
Finally, (iii) means that, for every non-zero $h\in\HH, \; E_{Xh}$ is a scalar spectral measure.
Hence all these localizations are equivalent, and so $\TT$ is quasianalytic.

\rightline{$\square$}      

A finer distinction can be made among Borel sets in the absolutely continuous (a.c.) case.
Let $m_n$ denote the normalized Lebesgue measure on $\T^n$.
(More precisely, we consider its restriction to $\B_n$.)
Let us assume that $\TT\in C_{*\cdot}$ is an \emph{a.c.\ commuting $n$-tuple}, which means  that $\m\ac m_n$, where $\m$ is a scalar spectral measure of $\UU$.
Taking the Radon--Nikodym derivative $0\le g_E= d\m/dm_n\in L^1(m_n)= L^1(\T^n)$, the Borel set $\o_a(\UU):= \{\zz\in \T^n: g_E(\zz)>0\}$ is called the \emph{ a.c.\ (global) residual set} of $\UU$.
For every $y\in\K, \; E_y\ac m_n$ also holds. 
Taking $0\le g_y = dE_y/dm_n\in L^1(\T^n)$, the Borel set $\o_a(\UU,y):= \{\zz\in\T^n: g_y(\zz)>0\}$ is called the \emph{ a.c.\ local residual set} of $\UU$ at $y$.
It is immediate that $E(\o_a(\UU,y))\K$ is the smallest spectral subspace, containing $y$.
It is clear also that $\o_a(\UU,y)\subset \o_a(\UU)$ for all $y\in\K$, and $\o_a(\UU,y)= \o_a(\UU)$ for some $y\in\K$ (see Proposition \ref{scalar-spm}).
We obtain that $\o_a(\UU)$ is the \emph{measurable support} of $E$ (with respect to $m_n$), that is  $E(\o)\ne 0$ holds exactly when $m_n(\o\cap \o_a(\UU))>0\; (\o\in\B_n)$.

For any $\o_1, \o_2\in \B_n$, the equivalence relation $\o_1 \overset{e}= \o_2$ holds by definition if $m_n(\o_1\triangle \o_2)=0$.
The relation $\o_1 \overset{e}\subset \o_2$, defined by $m_n(\o_1\setminus\o_2)=0$, is a partial ordering on $\B_n$, with the previous extended identity.

\begin{pro}
\label{ac-clattice}
The pair $(\B_n, \overset{e}\subset)$ is a complete lattice.
\end{pro}

\noindent
\textbf{Proof.}
Given any system $\{\o_j\}_{j\in J}$ in $\B_n$, let $\underline\O$ be the collection of those Borel sets $\o$, which satisfy the condition $\o \overset{e}\subset \o_j$ for all $j\in J$.
Setting $\underline\delta = \sup\{m_n(\o): \o\in\underline\O\}$, there exists a sequence $\{\o_{j_l}\}_{l\in\N}$ such that $\lim_{l\to\i} m_n(\o_{j_l})= \underline\d$.
It is easy to see that $\underline\o = \overset{e}\wedge\{\o_j\}_{j\in J}$ holds for the Borel set $\underline\o := \cup_{l=1}^\i \o_{j_l}$.
The existence of $\overset{e}\vee \{\o_j\}_{j\in J}$ can be shown similarly.

\rightline{$\square$}

For any $h\in\HH, \; \o_a(\TT,h):= \o_a(\UU,Xh)$ is called the \emph{a.c.\ local residual set} of $\TT$ at $h$.
It follows by the minimality of $(X,\UU)$ that the \emph{a.c. (global) residual set} $\o_a(\TT):= \overset{e}\vee\{\o_a(\TT,h) : h\in\HH\}$ of $\TT$ coincides with $\o_a(\UU)$.

The Borel set $\pi_a(\TT):= \overset{e}\wedge \{\o_a(\TT,h): 0\ne h\in\HH\}$ is called the \emph{ a.c. quasianalytic spectral set} of $\TT$.
The previous spectral invariants are determined up to sets of zero Lebesgue measure, and they are independent of the special choice of the unitary asymptote.
We note also that these definitions are compatible with those given in \cite{KerOT} for single a.c.\ polynomially bounded operators.

For any $\o_1, \o_2\in\B_n, \; \o_1 \subE\o_2$ holds if and only if $\o_1 \cap \o_a(\UU) \overset{e}\subset \o_2\cap \o_a(\UU)$.
Furthermore, $ \o_a(\TT,h) \eqE \o(\TT,h)$ is true for every $h\in\HH$.
It follows that $\o_a(\TT)\eqE \o(\TT) \eqE \T^n$ and $\pi_a(\TT) \eqE \pi(\TT)$.
Thus, the following statement can be derived from Theorem \ref{qa-coincidence}.

\begin{thm}
\label{ac-coincidence}
For the a.c.\ commuting $n$-tuple $\TT\in C_{*\cdot}$, the following conditions are equivalent:
\begin{itemize}
\item[\textup{(i)}] $\TT$ is quasianalytic;
\item[\textup{(ii)}] $\pi_a(\TT)= \o_a(\TT)$;
\item[\textup{(iii)}] $\o_a(\TT,h)= \o_a(\TT)$, for all non-zero $h\in\HH$.
\end{itemize}
\end{thm}

\medskip
\noindent
\textsc{L. K\'erchy}, Bolyai Institute, University of Szeged, Aradi v\'ertan\'uk tere 1, H-6720 Szeged, Hungary; \emph{e-mail}: kerchy@math.u-szeged.hu

\end{document}